\documentclass[12pt]{amsart}
\usepackage{amsmath, amssymb, amsthm}

\newtheorem{theorem}{Theorem}[section]
\newtheorem{definition}[theorem]{Definition}

\newtheorem{corollary}[theorem]{Corollary}
\newtheorem{lemma}[theorem]{Lemma}
\newtheorem{fact}[theorem]{Fact}
\newtheorem*{example}{Example}

\newtheorem*{sauer-shelah}{Sauer-Shelah Lemma}
\newtheorem*{volume-ratio}{Volume Ratio Theorem}

\def \proof {\noindent {\bf Proof.}\ \ }
\def \remark {\noindent {\bf Remark.}\ \ }
\def \remarks {\noindent {\bf Remarks.}\ \ }
\def \endproof {{\mbox{}\nolinebreak\hfill\rule{2mm}{2mm}\par\medbreak}}

\DeclareMathOperator{\codim}{codim} 
\DeclareMathOperator{\rank} {rank}

\def \R {\mathbb{R}}

\def \Z {\mathbb{Z}}
\def \E {\mathbb{E}}

\def \CC {\mathcal{C}}
\def \EE {\mathcal{E}}

\def \e {\varepsilon}

\def \d {\delta}

\def \l {\lambda}
\def \L {\Lambda}

\def \< {\langle}
\def \> {\rangle}

\def \diam {{\rm diam}}

\def \conv {{\rm conv}}
\def \cconv {{\rm cconv}}

\def \vol {{\rm vol}}

\begin{document}

\title {Integer cells in convex sets}

\author {R. Vershynin}

\address{Department of Mathematics, University of California, Davis,
  1 Shields Ave, Davis, CA 95616, U.S.A.}
\email{vershynin@math.ucdavis.edu}

\thanks{partially supported by the New Faculty Research Grant 
  of the University of California-Davis and by the NSF Grant 0401032.}



\date{November 4, 2004}

\maketitle

\begin{abstract}  
  Every convex body $K$ in $\R^n$ has a coordinate projection $PK$ that 
  contains at least $\vol(\frac{1}{6} K)$ cells of the integer lattice $P\Z^n$, 
  provided this volume is at least one.
  Our proof of this counterpart of Minkowski's theorem 
  is based on an extension of the combinatorial density 
  theorem of Sauer, Shelah and Vapnik-Chervonenkis to $\Z^n$. 
  This leads to a new approach to sections of convex bodies.
  In particular, fundamental results of the asymptotic convex geometry
  such as the Volume Ratio Theorem and Milman's duality of the diameters 
  admit natural versions for coordinate sections. 
\end{abstract}


\section{Introduction}

Minkowski's Theorem, a central result in the geometry of numbers, 
states that if $K$ is a convex and symmetric set in $\R^n$, 
then $\vol(K) > 2^n$ implies that $K$ contains a
nonzero integer point. More generally, $K$ contains at least
$\vol(\frac12K)$ integer points. The main result of the present paper
is a similar estimate on the number of {\em integer cells}, 
the unit cells of the integer lattice $\Z^n$, contained in a convex body.

Clearly, the largeness of the volume of $K$ does not imply the existence 
of any integer cells in $K$; a thin horizontal pancake is an example.
The obstacle in the pancake $K$ is caused by only one coordinate in which 
$K$ is flat; after eliminating it (by projecting $K$ onto the remaining
ones) the projection $PK$ will have many integer cells 
of the lattice $P\Z^n$. This observation turns out to be a general 
principle. 

\begin{theorem}                              \label{combin volume}
  Let $K$ be a convex set in $\R^n$.
  Then there exists a coordinate projection $P$
  such that $PK$ contains at least $\vol(\frac{1}{6} K)$
  cells of the integer lattice $P\Z^n$, 
  provided this volume is at least one.
\end{theorem}
\noindent (A coordinate projection is the orthogonal projection in $\R^n$
onto $\R^I$ for some nonempty subset $I \subset \{1,\ldots,n\}$.)

\medskip
\paragraph{Combinatorics: Sauer-Shelah-type results.}
Theorem \ref{combin volume} is a consequence of an extension to $\Z^n$
of the famous result due to Vapnik-Chervonenkis, Sauer, 
Perles and Shelah, commonly known as Sauer-Shelah lemma,
see e.g. \cite[\S 17]{B}.

\begin{sauer-shelah}
  If $A \subset \{0,1\}^n$ has cardinality 
  $\# A > \binom{n}{0} +  \binom{n}{1} + \ldots +  \binom{n}{d}$,
  then there exists a coordinate projection $P$ of rank larger than
  $d$ and such that $PA = P\{0,1\}^n$.
\end{sauer-shelah}
\noindent This result is used in a variety of areas 
ranging from logics to theoretical computer science
to functional analysis \cite{M 3}.
In order to bring Sauer-Shelah Lemma to geometry, 
we will need first to generalize it to sets $A \subset \Z^n$. 
An {\em integer box} is a subset of $\Z^I$ of the form 
$\prod_{i \in I} \{a_i,b_i\}$ with $a_i \ne b_i$.
 
\begin{theorem}                                     \label{SSL Zn intro}
  If $A \subset \Z^n$, then 
  $$
  \# A  \le  1 + \sum_P 
            \# \big( \text{integer boxes in $PA$} \big),
  $$
  where the sum is over all coordinate projections $P$.
\end{theorem}
\noindent If $A \subset \{0,1\}^n$, then every $PA$ in the sum above
may contain only one integer box $P\{0,1\}^n$ if any, hence
\begin{equation}                                    \label{more than SSL intro}
  \# A  \le  1 +  \# \big( P \text{ for which } PA = P\{0,1\}^n \big).
\end{equation} 
Estimate \eqref{more than SSL intro} is due to A.Pajor \cite{Pa}.
Since the right hand side of \eqref{more than SSL intro} is bounded by 
$\binom{n}{0} +  \binom{n}{1} + \ldots +  \binom{n}{d}$,
where $d$ is the maximal rank of $P$ for which $PA = P\{0,1\}^n$,
\eqref{more than SSL intro} immediately implies Sauer-Shelah Lemma.

In a similar way, Theorem \ref{SSL Zn intro} implies a recent 
generalization of Sauer-Shelah lemma in terms of Natarajan dimension, 
due to Haussler and Long \cite{HL}. 
In their result, $A$ has to be bounded
by some paralelepiped; we do not impose any boundedness restrictions
(see Corollary \ref{HL}).

Most importantly, Theorem \ref{SSL Zn intro} admits a version 
for integer {\em cells} instead of integer boxes. If $A \subset \R^n$ 
is convex, then 
$$
\# A  \le  1 + \sum_P 
        \# \big( \text{integer cells in $PA$} \big).
$$
This quickly leads to Theorem \ref{combin volume}.
This version also implies a generalization of Sauer-Shelah lemma 
from \cite{HL} in terms of the {\em combinatorial dimension}, 
which is an important concept originated in the statistical 
learning theory and and which became widely useful in many areas, 
see \cite{ABCH}, \cite{CH}, \cite{Ta}, \cite{MV}.
These results will be discussed in detail in Section \ref{s:combinatorics}. 
The proof of Theorem~\ref{SSL Zn intro}
relies on the combinatorics developed in \cite{MV} and \cite{RV}.

\medskip
\paragraph{Convex geometry: coordinate sections of convex bodies.}
Theorem~\ref{SSL Zn intro} leads to a new approach to 
coordinate sections of convex bodies.

The problem of finding nice coordinate sections of a symmetric 
convex body $K$ in $\R^n$ has been extensively studied in geometric 
functional analysis. 
It is connected in particular with important 
applications in harmonic analysis, where the system of characters 
defines a natural coordinate structure. 
The $\L_p$-problem, which was solved by J.~Bourgain \cite{Bou},
is an exemplary problem on finding nice coordinate sections,
as explained by an alternative and more general solution 
(via the majorizing measures) given by M.~Talagrand \cite{Ta 95}.
It is generally  extremely difficult to find a nice coordinate section 
even when the existence of nice generic sections (usually randomly 
chosen from the Grassmanian) is well known, see e.g. \cite{Ta}, 
\cite{MV}, \cite{RV}.

The method of the present paper
allows one to prove natural versions of a few classical results for 
coordinate sections. 
Since the number of integer cells in a set $K$ is bounded by its 
volume, we have in Theorem \ref{combin volume} that 
\begin{equation}                                          \label{PK K}
  \text{$PK$ contains an integer cell and 
        $|PK| \ge |{\textstyle \frac{1}{6}} K|$.}
\end{equation}
(we write $|PK| = \vol(PK)$ for the volume in $P\R^n$).
This often enables one to conclude {\em a posteriori} that $P$
has large rank, as \eqref{PK K} typically fails for 
all projections of small ranks.

If $K$ is symmetric and an integer $m <n$ is fixed, 
then using \eqref{PK K} for $a^{-1} K$ 
with an appropriate $a > 0$, we obtain
$a^{-m} |PK| \ge  a^{-n} |\frac{1}{6} K|$ 
for some coordinate projection $P$ of rank $m$.
Moreover,  $P(a^{-1}K)$ contains a unit coordinate cube, 
so solving for $a$ we conclude that
\begin{equation}                                          \label{ratio intro}
  \text{$PK$ contains a coordinate cube of side } 
  \left( \frac{|c K|}{|P K|} \right)^{\frac{1}{n-m}}. 
\end{equation}
where $C, c, c_1, \ldots$ denote positive absolute constants 
(here $c = 1/6$). 

This leads to a ``coordinate'' version of the classical
Volume Ratio Theorem. 
This theorem is a  remarkable phenomenon originated
in the work of B.~Kashin related to approximation theory \cite{K}, 
developed by S.~Szarek into a general method \cite{Sz} and carried over
to all convex bodies by S.~Szarek and N.~Tomczak-Jaegermann 
(\cite{STo}, see \cite[\S 6]{Pi}). 
The unit ball of $L_p^n$ ($1 \le p \le \infty$) is denoted
by $B_p^n$, i.e. for $p < \infty$
$$
x \in B_p^n \ \ \text{iff} \ \ |x(1)|^p + \cdots + |x(n)|^p \le n
$$
and $x \in B_\infty^n$ if $\max_i |x(i)| \le 1$.
We choose to work with $L_p^n$ spaces rather than $\ell_p^n$ spaces here
because their unit balls have approximately unit volume: 
$c_1 \le |B_p^n|^{1/n} \le c_2$ for all $1 \le p \le \infty$.

\begin{volume-ratio}(Szarek, Tomczak-Jaegermann).
  Let $K$ be a convex symmetric body in $\R^n$ which contains $B_2^n$. 
  Then for every integer $0 < k < n$ there exists 
  a subspace $E$ of codimension $k$ and such that 
  \begin{equation}                                  \label{vr}
  K \cap E  \subseteq |C K|^{1/k} B_2^n.
  \end{equation}
\end{volume-ratio}
\noindent In fact, the subspace $E$ can be taken at random from the Grassmanian.

\medskip To obtain a coordinate version of the Volume Ratio Theorem, 
we can not just claim that \eqref{vr} holds for some 
coordinate subspace $E = \R^I$: the octahedron $K = B_1^n$ 
forms an obstacle. However it turns out that the octahedron 
is the only obstacle, so our claim becomes true if one 
replaces the Euclidean ball $B_2^n$ in \eqref{vr} by 
its circumscribed octahedron $B_1^n$. This seems to be a general 
phenomenon when one passes from arbitrary to coordinate
sections, see \cite{RV}.

\begin{theorem}                               \label{cvr cube inside intro} 
  Let $K$ be a convex symmetric body in $\R^n$ which contains $B_\infty^n$. 
  Then for every integer $0 < k < n$ there exists 
  a coordinate subspace $E$ of codimension $k$ and such that 
  $$
  K \cap E  \subseteq |C K|^{1/k} B_1^n.
  $$
\end{theorem}
\noindent This theorem follows from \eqref{ratio intro} by duality
(Santalo and the reverse Santalo inequalities, the latter due to 
Bourgain and Milman).

\medskip
\remarks {\bf 1.} The assumption $B_\infty^n \subset K$ of 
Theorem \ref{cvr cube inside intro} is weaker than 
the assumption $B_2^n \subset K$ of the Volume Ratio Theorem.
In fact, this assumption can be completely eliminated if one 
replaces $|CK|^{1/k}$ by the quantity
$$
A_k(K) = \max \left( \frac{|CK|}{|K \cap E|} \right)^{1/\codim E}
$$
where the maximum is over the coordinate subspaces $E$, $\codim E \ge k$.
Clearly, $A_k(K) \le |CK|^{1/k}$ if $K$ contains $B_\infty^n$.
We will discuss this ``Coordinate Volume Ratio Theorem''
as well as the quantity $A_k(K)$ in more detail in Section~\ref{s:cvr}.

{\bf 2.} The right dependence on $k/n$ in the Volume Ratio Theorem and 
in Theorem \ref{cvr cube inside intro} is a delicate problem.
$|CK|^{1/k} = C^{n/k} |K|^{1/k}$, and while the factor $|K|^{1/k}$
is sharp (which is easily seen for ellipsoids or parallelepipeds),
the exponential factor $C^{n/k}$ is {\em not}. We will improve
it (in the dual form) to a linear factor $Cn/k$ in 
Section~\ref{s:Lp}. 

\medskip
Another example of applications of Theorem \ref{combin volume}
is a coordinate version of Milman's duality of diameters 
of sections. For a symmetric convex body $K$ in $\R^n$, let 
$$
b_k(K) = \min \diam(K \cap E_k),
$$
where the minimum is over all $k$-dimensional subspaces $E_k$. 
Then for every $\e > 0$ 
and for any two positive integers $k$ and $m$ satisfying 
$k + m \le (1-\e)n-C$ one has 
\begin{equation}                            \label{milman intro}
b_k(K) \, b_m(K^\circ) \le C/\e.
\end{equation}
(in fact, this holds for random subspaces $E_k$ in the Grassmanian)
\cite{Mi 1}, \cite{Mi 2}.
This phenomenon reflects deep linear duality relations 
and provides a key tool in understanding the ``global''
duality in asymptotic convex geometry, see \cite{Mi 2}, 
\cite{Mi 98}.

To establish a version of this result for coordinate subspaces $E_k$,
we have (as before) to change the metric that defines the diameter 
to that given by the octahedron circumscribed around the unit
Euclidean ball (rather than the Euclidean ball itself).
Then for the new diameter $\diam_1$ we let 
$$
r_k(K) = \min \diam_1(K \cap E_k),
$$
where the minimum is over all $k$-dimensional {\em coordinate} 
subspaces $E_k$. 
In other words, the inequality $r_k(K) \le 2r$ holds 
iff one can find a $k$-element set $I$ so that one has
$\sum_{i \in I} |x(i)|  \le  r \sqrt{n}$ for all $x \in K$.

\begin{theorem}[Duality for diameters of coordinate sections] 
       \label{coordinate}
  Let $K$ be a symmetric convex body in $\R^n$.
  For any $\e > 0$ and for any two positive integers $k$ and $m$ 
  satisfying $k + m \le (1-\e)n$
  one has 
  $$
  r_k(K) \, r_m(K^\circ) \le C^{1/\e}.
  $$
\end{theorem}
\noindent In particular, there exists a subset of coordinates $I$ of size,
say, $\lceil n/3 \rceil$
such that the absolute values of the coordinates in $I$ sum to at most 
$C \sqrt{n}$ either for all vectors in $K$ or for all vectors
in $K^\circ$.

\medskip
\remark In most of the results of this paper, 
the convexity of $K$ can be relaxed to a weaker  
{\em coordinate convexity}, see e.g. \cite{M 2}.

\qquad

ACKNOWLEDGEMENTS. The author is grateful to M. Rudelson for his interest
and valuable comments. This project started when the author was 
at the Pacific Institute for Mathematical Studies and the 
University of Alberta. He thanks these institutions 
and especially N. Tomczak-Jaegermann for support.

\section{Sauer-Shelah Lemma in $\Z^n$}              \label{s:combinatorics}

In 1971-72, Vapnik and Chervonenkis \cite{VC}, Sauer \cite{Sa} and 
Perles and Shelah \cite{Sh} independently proved the following well
known result, which has found applications in a variety of areas
ranging from logics to probability to computer science.

\begin{theorem}[Sauer-Shelah Lemma]
  If $A \subset \{0,1\}^n$ has cardinality 
  $\# A > \binom{n}{0} +  \binom{n}{1} + \ldots +  \binom{n}{d}$,
  then there exists a coordinate projection $P$ of rank larger than
  $d$ and such that 
  \begin{equation}                                   \label{SSL PA}
    PA = P\{0,1\}^n.
  \end{equation}
\end{theorem}
\noindent A short proof of Sauer-Shelah Lemma can be found e.g. in \cite[\S 17]{B};
for numerous variants of the Lemma see the bibliography in \cite{HL}
as well as \cite{Ale}, \cite{ST 1}, \cite{ST 2}.

\medskip
To bring Sauer-Shelah Lemma to geometry, 
we will have to generalize it to sets $A \subset \Z^n$.
The case when such $A$ is bounded by a parallelopiped, 
i.e. $A \subset \prod_{i=1}^n \{0,\ldots,N_i\}$, 
is well understood by now, see \cite{KM}, \cite{Alo}, \cite{St}, 
\cite{HL}. In this section we will prove a generalization 
of Sauer-Shelah Lemma to $A \subset \Z^n$ independent of any boundedness
assumptions.

We start with a simpler result. 
An {\em integer box} is a subset of $\Z^n$ of the form 
$\{a_1, b_1\} \times \cdots \times \{a_n, b_n\}$ with 
$a_i \ne b_i$ $\forall i$. Similarly one defines 
integer boxes in $\Z^I$, where $I \subseteq \{1,\ldots,n\}$.

\begin{theorem}                                     \label{SSL Zn}
  If $A \subset \Z^n$, then 
  \begin{equation}                                  \label{SSL A}
    \# A  \le  1 + \sum_P 
            \# \big( \text{integer boxes in $PA$} \big),
  \end{equation}
  where the sum is over all coordinate projections $P$ in $\R^n$.
\end{theorem}

\remark
Let $A \subset \{0,1\}^n$. 
Since the only lattice box that can be contained in $PA$ is $P\{0,1\}^n$, 
Theorem~\ref{SSL Zn} implies that
\begin{equation}                                    \label{more than SSL}
  \# A  \le  1 +  \# \big( P \text{ for which } PA = P\{0,1\}^n \big).
\end{equation} 
This estimate is due to A.Pajor \cite[Theorem 1.4]{Pa}.
Note that this quantity is bounded by 
$\binom{n}{0} +  \binom{n}{1} + \ldots +  \binom{n}{d}$,
where $d$ is the maximal rank of $P$ for which $PA = P\{0,1\}^n$.
This immediately implies Sauer-Shelah Lemma.

\medskip
The result that we really need for geometric applications is 
Theorem \ref{SSL Zn} for {\em integer cells}, 
which are integer boxes whose all sides equal $1$.
Although the number of integer cells in a convex body can in 
principle be estimated through the number of integer boxes, 
the dependence will not be linear -- a cube $[0,M]^n$ contains
$M^n$ integer cells and $(\frac{1}{2}M(M+1))^n$ integer boxes.
To obtain Theorem \ref{SSL Zn} for integer cells, we will have 
to prove a more accurate extension of Sauer-Shelah Lemma to $\Z^n$. 

The crucial in our discussion will be the notion of 
{\em coordinate convexity} (see e.g. \cite{M 2}),
which is weaker than that of convexity.

\begin{definition}
  Let $K$ be a set in $\R^n$. The {\em coordinate convex hull} of $K$
  consists of the points $x \in \R^n$ such that for every choice of signs
  $\theta \in \{-1,1\}^n$ one can find $y \in K$ such that
  \begin{align*}
    y(i) \ge x(i) &\ \ \ \text{if $\theta(i) = 1$,} \\
    y(i) \le x(i) &\ \ \ \text{if $\theta(i) = -1$.}
  \end{align*}
  $K$ is called {\em coordinate convex} if it coincides with its
  coordinate convex hull.
\end{definition}
By changing $\R^n$ to $\Z^n$ the coordinate convexity can also 
be defined for subsets of $\Z^n$. Also, changing $\R^n$ to $\R^I$ and 
$\{-1,1\}^n$ to $\{-1,1\}^I$, the coordinate convexity is defined 
for subsets of $\R^I$ (and similarly for $\Z^I$), 
where $I \subset \{1,\ldots,n\}$.

One obtains a general convex body in $\R^n$ by cutting off half-spaces. 
Similarly, a general coordinate convex body in $\R^n$ is obtained by 
cutting off octants, i.e. translates of the sets 
$\theta \cdot \R_+^n$ with $\theta \in \{-1,1\}^n$.
Clearly, every convex set is coordinate convex;
the converse is not true, as the cross shows
$\{ (x,y) \; | \; x = 0 \text{ or } y = 0 \}$ in $\R^2$.

The central combinatorial result of this section is the following theorem
which we will prove after some comments.

\begin{theorem}                                   \label{lattice vol} 
  For every $A \subset \Z^n$,
  \begin{equation}                                \label{lattice vol eqn}
  \# A  \le  1 + \sum_P  
         \# \big( \text{integer cells in $\cconv PA$} \big),
  \end{equation} 
  where the sum is over all coordinate projections $P$.
\end{theorem}

\paragraph{The combinatorial dimension and Sauer-Shelah type results.}
Like Theorem \ref{SSL Zn}, Theorem \ref{lattice vol} also contains 
Sauer-Shelah Lemma: 
every subset $A \subset \{0,1\}^n$ is coordinate convex, 
and the only lattice box that can be contained in $PA$ is $P\{0,1\}^n$, 
which implies \eqref{more than SSL} and hence Sauer-Shelah lemma.

To see a relation of Theorem \ref{lattice vol} to later generalizations
of Sauer-Shelah lemma, let us recall an important concept of the 
combinatorial dimension, which originates in the statistical learning 
theory and which became useful in convex geometry, combinatorics 
and analysis, see \cite{ABCH}, \cite{CH}, \cite{Ta}, \cite{MV}, \cite{RV}.

\begin{definition}
  The {\em combinatorial dimension} $v(A)$ of a set $A \subset \R^n$
  is the maximal rank of a coordinate projection $P$ such that 
  $\cconv(PA)$ contains some translate of the unit cube $P\{0,1\}^n$.
  
  For $t > 0$, the scale-sensitive version of the combinatorial dimension
  is defined as $v(A,t) = v(t^{-1}A)$.
\end{definition}
Equivalently, a subset $I \subset \{1,\ldots, n\}$ is called $t$-shattered 
by $A$ if there exists an $h \in \R^n$ such that, given any partition 
$I = I^- \cup I^+$, one can find an $x \in A$ such that 
$x(i) \le h(i)$ if $i \in I^-$ and $x(i) \ge h(i)+t$ if $i \in I^+$.
The combinatorial dimension $v(A,t)$ is the maximal cardinality of a 
subset $t$-shattered by $A$.

\medskip
A few words on the history of the concept of the combinatorial dimension.
For sets $A \subset \{0,1\}^n$, the combinatorial dimension $v(A)$ 
is the classical {\em Vapnik-Chernovenkis dimension}; see \cite{M 3} 
for a nice introduction to this important concept. 
For sets $A \subset \Z^n$, the notion of the combinatorial dimension $v(A)$ 
goes back to 1982-83, when Pajor used it for origin symmetric classes 
in view of applications to the local theory of Banach spaces \cite{Pa seminar}. 
He proved early versions of Sauer-Shelah Lemma for sets $A \subset \{0,\ldots,p\}^n$ 
(see \cite{Pa seminar}, \cite[Lemma 4.9]{Pa}).
Pollard gave an explicit definition of $v(A)$ in his 1984 book on 
stochastic processes \cite{Po}.
Haussler also discussed this concept in his 1989 work in learning theory 
(\cite{Ha PAC}, see \cite{HL} and the references therein).  

For convex and origin symmetric sets $A \subset \R^n$, the combinatorial 
dimension $v(A,t)$ is easily seen to coincide with the maximal rank of the
coordinate projection $PA$ of $A$ that contains the centered coordinate 
cube of side $t$. In view of this straightforward connection to convex geometry 
and thus to the local theory of Banach spaces, the combinatorial dimension 
was a central quantity in the 1982-83 works of Pajor (\cite {Pa volumes}, 
see Chapter IV of \cite{Pa}).
Connections of $v(A,t)$ to Gaussian processes and further applications to 
Banach space theory were established in the far reaching 1992 paper of 
M.Talagrand (\cite{Ta 92}, see also \cite{Ta}). 
The quantity $v(A,t)$ was formally defined in 1994 by Kearns and Schapire 
for general sets $A$ in their paper in learning theory \cite{KS}.

Since its invention, the combinatorial dimension turned out to be very effective 
in measuring the complexity of a set $A$ in combinatorics, 
functional analysis, statistical learning theory, the theory of empirical processes, 
discrete and convex geometry
(see \cite{ABCH}, \cite{Ta}, \cite{MV}, \cite{RV}).
Alternative names for the combinatorial dimension used in the literature on 
combinatorics and statistical learing theory are:
Pollard dimension and pseudo dimension for $v(A)$, 
shattering and fat-shattering dimension for $v(A,t)$, see \cite{HL} and 
\cite{M fewnotes}.

Similarly, {\em Natarajan dimension} $n(A)$ of a set $A \subset \Z^n$
is the maximal rank of a coordinate projection $P$ such that 
$PA$ contains an integer box (see \cite{HL}).

\medskip
Theorems \ref{SSL Zn} and \ref{lattice vol} easily imply two 
recent results of Haussler and Long \cite{HL} on 
the combinatorial and Natarajan dimensions, 
which are in turn generalizations of Sauer-Shelah lemma. 

\begin{corollary}[Haussler, Long \cite{HL}]				\label{HL}
  Let $A \subset \prod_{i=1}^n \{0, \ldots, N_i\}$. Then
  
  (i) We have
  $$
  |A|  \le  \sum_{\# I \le v(A)} \prod_{i \in I} N_i,
  $$
  where the sum is over all subsets $I \subset \{1,\ldots,n\}$ 
  of cardinality at most $v(A)$ (we include $I = \emptyset$ and assign
  to it the summand equal to $1$).
  
  (ii) In particular, if $A \subset \{0,\ldots, N\}^n$ then 
  $$
  |A|  \le  \sum_{i=0}^{v(A)} \binom{n}{i} N^i.
  $$
  
  (iii) We have 
  $$
  |A|  \le  \sum_{\# I \le n(A)} 
                   \prod_{i \in I} \binom{N_i + 1}{2}.
  $$
  where the sum is over all subsets $I \subset \{1,\ldots,n\}$ 
  of cardinality at most $n(A)$ (we include $I = \emptyset$ and assig
  to it the summand equal to $1$). 
\end{corollary}
  
\proof
For (i), apply Theorem \ref{lattice vol}. All the summands in 
\eqref{lattice vol eqn} that correspond to $\rank P > v(A)$ 
vanish by the definition of the combinatorial dimension.
Each of the non-vanishing summands is bounded by the number 
of integer cells in $\cconv PA \subset P(\prod_{i=1}^n \{0, \ldots, N_i\})$.
This establishes (i) and thus (ii).

Repeating this for (iii), we only have to note that the number of 
integer boxes in 
$P(\prod_{i=1}^n \{0, \ldots, N_i\}) 
= \{0\} \times \prod_{i \in I} \{0, \ldots, N_i\}$ is at 
most $\binom{N_i + 1}{2}$.
\endproof
 
\medskip \remark All the statements in Corollary \ref{HL}
reduce to Sauer-Shelah lemma if $N_i = 1$ $\forall i$.

\paragraph{The proof.}
Here we prove Theorem \ref{lattice vol}.
Define the {\em cell content of $A$} as  
$$
\Sigma(A) = \sum_P 
      \# \big( \text{integer cells in $\cconv PA$} \big),
$$
where we include in the counting one $0$-dimensional projection $P$ 
(onto $\R^\emptyset$),
for which the summand is set to be $1$ if $A$ is nonempty and $0$ otherwise.
This definition appears in \cite{RV}. 
We partition $A$ into sets $A_k$, $k \in \Z$, defined as 
$$
A_k = \{x \in A :\; x(1) = k\}.
$$

\begin{lemma}                       \label{induction}
  For every $A \subset \Z^n$, 
  $$
  \Sigma(A)  \ge  \sum_{k \in \Z} \Sigma(A_k).
  $$
\end{lemma}

\proof
A cell $\CC$ in $\R^I$, $I \subset \{1,\ldots,n\}$, will be considered
as an ordered pair $(\CC,I)$. This also applies to the trivial cell 
$(0,\emptyset)$ which we will include in the counting throughout
this argument.
The coordinate projection onto $\R^I$ will be denoted by $P_I$.

We say that $A$ {\em has a cell} $(\CC,I)$ if 
$\CC \subset \cconv P_I B$.
The lemma states that $A$ has at least as many cells as all the sets $A_k$
have in total.

If $A_k$ has a cell $(\CC,I)$ then $A$ has it, too.
Assume that $N>1$ sets among $A_k$ have a nontrivial cell $(\CC,I)$.
Since the first coordinate of any point in such a set $A_k$ equals $k$,
one necessarily has $1 \not\in I$. Then 
$P_{\{1\} \cup I} A_k = \{k\} \times P_I A_k$, where the factor $\{k\}$
means of course the first coordinate. Hence
\begin{align*}
\{k\} \times \CC 
&\subset \{k\} \times \cconv(P_I A_k)
= \cconv (\{k\} \times P_I A_k) \\
&= \cconv P_{\{1\} \cup I} A_k 
\subset \cconv P_{\{1\} \cup I} A.
\end{align*}
Therefore the set $\cconv P_{\{1\} \cup I} A$ contains the 
integer box $\{k_1,k_2\} \times \CC$, where $k_1$ is the 
minimal $k$ and $k_2$ is the maximal $k$ for the $N$ sets $A_k$.
Then $\cconv P_{\{1\} \cup I} A$ must also contain 
$\cconv(\{k_1,k_2\} \times \CC) \supset [k_1, k_2] \times \CC$
which in turn contains at least $k_2 - k_1 \ge N-1$ integer cells
of the form $\{a,a+1\} \times \CC$.
Hence, in addition to one cell $\CC$, the set $A$ has at least
$N-1$ cells of the form 
\begin{equation}                    \label{new cell}
  (\{a,a+1\} \times \CC, \{1\} \cup I).
\end{equation}
Since the first coordinate of all points in any fixed $A_k$ is 
the same, none of $A_k$ may have a cell of the form \eqref{new cell}.
Note also that the argument above works also for the trivial cell.

This shows that there exists an injective mapping from the set of
the cells that at least one $A_k$ has into the set of the cells that 
$A$ has. The lemma is proved.
\endproof

\medskip
\noindent {\bf Proof of Theorem \ref{lattice vol}.}
It is enough to show that for every $A \subset \Z^n$ 
$$
\# A \le \Sigma(A).
$$
This is proved using Lemma \ref{induction} by induction 
on the dimension $n$.

The claim is trivially true for $n=0$ (in fact also for $n=1$).
Assume it is true for some $n \ge 0$.
Apply Lemma \ref{induction} and note that each $A_k$ is a translate
of a subset in $\Z^{n-1}$. We have
$$
\Sigma(A)  \ge  \sum_{k \in \Z} \Sigma(A_k)
\ge \sum_{k \in \Z} \# A_k 
= \# A
$$
(here we used the induction hypothesis for each $A_k$).
This completes the proof.
\endproof

\paragraph{Volume and lattice cells}
Now we head to Theorem \ref{combin volume}.

\begin{corollary}               \label{vol cell content}
  Let $K$ be set in $\R^n$.
  Then 
  $$
  |{\textstyle\frac{1}{2}} K|  
  \le  1 + \sum_P  
         \# \big( \text{integer cells in $\cconv PK$} \big),
  $$ 
  where the sum is over all coordinate projections $P$.
\end{corollary}

\noindent For the proof we need a simple fact:

\begin{lemma}                   \label{translate}
  For every set $K$ in $\R^n$ and every $x \in \R^n$, 
  $$
  \# \big( \text{integer cells in $x+K$}  \big)
  \le \# \big( \text{integer cells in $2K$} \big). 
  $$
\end{lemma}

\proof
The proof reduces to the observation that every translate of 
the cube $[0,2]^n$ by a vector in $\R^n$ contains an integer cell.
This in turn is easily seen by reducing to the one-dimensional case.
\endproof

\noindent{\bf Proof of Corollary \ref{vol cell content}.}
Let $x$ be a random vector uniformly distributed in $[0,1]^n$, 
and let $A_x = (x+K) \cap \Z^n$.
Then $\E \# A_x = |K|$.
By Theorem \ref{lattice vol}, 
\begin{equation}                \label{average content}
|K| \le 1 + \E \sum_P  
         \# \big( \text{integer cells in $\cconv PA_x$} \big),
\end{equation}
while 
\begin{equation}                \label{cconv conv}
\cconv PA_x  \subset \cconv P(x+K) = Px + \cconv PK.
\end{equation}
By this and Lemma \ref{translate},
$$
\# \big( \text{integer cells in $\cconv PA_x$} \big) 
\le \# \big( \text{integer cells in $\cconv P(2K)$} \big). 
$$
Thus by \eqref{average content}
$$
|K| \le 1 + \sum_P  
         \# \big( \text{integer cells in $\cconv P(2K)$} \big).
$$
This proves the corollary.
\endproof

\medskip \remark The proof of Theorem \ref{SSL Zn} is very similar 
and in fact is simpler than the argument above. One looks at 
$\Sigma(A) = \sum_P \# \big( \text{integer boxes in $PA$} \big)$
and repeats the proof without worrying about coordinate convexity.

\qquad

Now we can prove the main geometric result of this section.

\begin{theorem}                 \label{counting cells}
  Let $K$ be a set in $\R^n$.
  Then there exists a coordinate projection $P$ in $\R^n$
  such that $\cconv PK$ contains at least $|\frac{1}{4}K| - 2^{-n}$ 
  integer cells.
\end{theorem}

\proof
By Corollary \ref{vol cell content}, 
$$
|{\textstyle\frac{1}{2}} K|  
  \le  1 + (2^n-1) \max_P
         \# \big( \text{integer cells in $\cconv PK$} \big).
$$
Hence $\max_P \# \big( \text{integer cells in $\cconv PK$} \big)
\ge |\frac{1}{4}K| - 2^{-n}$.
\endproof

\medskip Note that $|\frac{1}{4}K| - 2^{-n} \ge |\frac{1}{6}K|$
if $|\frac{1}{6}K| \ge 1$. This implies Theorem \ref{combin volume}.

\section{The Coordinate Volume Ratio Theorem}           \label{s:cvr}

Let $K$ be a set in $\R^k$. For $0 < k < n$, define
$$
A_k(K) = \max \left( \frac{|CK|}{|K \cap E|} \right)^{1/\codim E}
$$
where the maximum is over the coordinate subspaces $E$, $\codim E \ge k$,
and $C>0$ is an absolute constant whose value will be discussed later.

\begin{theorem}[Coordinate Volume Ratio Theorem]        \label{cvr}
  Let $K$ be a convex symmetric set in $R^n$. 
  Then for every integer $0 < k < n$ there exists a coordinate 
  section $E$, $\codim E = k$, such that  
  $$
  K \cap E \subset A_k(K) \, B_1^n.
  $$
\end{theorem}
\noindent The proof relies on the extension on Sauer-Shelah Lemma in $\Z^n$
from the previous section and on the duality for the volume, 
which is Santalo and the reverse Santalo
inequalities (the latter due to J.Bourgain and V.Milman). 
We will prove Theorem \ref{cvr} in the end of this section.

\paragraph{1.} 
In the important case when $K$ contains the unit cube, we have 
$A_k(K) \le |CK|^{1/k}$. This implies:

\begin{corollary}                   \label{cvr cube inside} 
  Let $K$ be a convex body in $\R^n$ which contains 
  the unit cube $B_\infty^n$. 
  Then for every integer $0 < k < n$ there exists 
  a coordinate subspace $E$ of codimension $k$ and such that 
  $$
  K \cap E  \subseteq |C K|^{1/k} \, B_1^n.
  $$
\end{corollary}

The assumptions of this corollary are weaker than those of the 
classical Volume Ratio Theorem stated in the introduction, 
because the cube $B_\infty^n$ is inscribed into the Euclidean ball $B_2^n$.
The conclusion of Corollary \ref{cvr cube inside} 
is that some coordinate section $K \cap E$ is bounded by the 
octahedron $B_1^n$, which is circumscribed around the Euclidean ball $B_2^n$.
No stronger conclusion is for a {\em coordinate} section is possible:
$K = B_1^n$ itself is an obstacle.

Nevertheless, by a result of Kashin (\cite{K}, see a sharper estimate
in Garnaev-Gluskin \cite{GG}) a random section of $B_1^n$ in the 
Grassmanian $G_{n,k}$ with $k = \lceil n/2 \rceil$ is equivalent 
to the Euclidean ball $B_2^k$. Thus a random (no longer coordinate) 
section of $K \cap E$ of dimension, say, $\frac{1}{2} \dim(K \cap E)$
will already be a subset of $|CK|^{1/k} B_2^n$.
This shows that Corollary \ref{cvr cube inside} is close in nature
to the classical Volume Ratio Theorem. It gives {\em coordinate}
subspaces without sacrificing too much of the power of the Volume 
Ratio Theorem.\footnote{Even though in the Coordinate Volume Ratio Theorem
the coordinate section can not be random in general, 
a very recent work of Giannopoulos, Milman and Tsolomitis \cite{GMT}
and of the author \cite{V}
suggests that one can automatically regain randomness 
of a bounded section in the Grassmanian if one only knows the 
{\em existence} of a bounded section in the Grassmanian.}  

In the next section we will prove a (dual) result even sharper than 
Corollary \ref{cvr cube inside}. 

\paragraph{2.} 
The quantity $A_k(K)$ is best illustrated on the example of classical bodies.
If $K$ is the parallelopiped $\prod_{i=1}^n [-a_i,a_i]$ with 
semiaxes $a_1 \ge a_2 \ge \cdots \ge a_n >0$, then 
\begin{equation}                    \label{Ak axes}
  A_k(K) = (2C)^{n/k} \Big( \prod_{i=1}^k a_i \Big)^{1/k},
\end{equation}
a quantity proportional to the geometric mean of the largest
$k$ semiaxes. The same holds if $K$ is the ellipsoid with the 
coordinate nonincreasing semiaxes $a_i \sqrt{n}$, i.e. 
$x \in K$ iff $\sum_{i=1}^n x(i)^2/a_i^2 \le n$.
This is clearly better than
$$
|CK|^{1/k} = (2C)^{n/k} \Big( \prod_{i=1}^n a_i \Big)^{1/k},
$$
which appears in the classical Volume Ratio Theorem (note that 
the inclusion $B_2^n \subset K$ implies in the ellipsoidal example 
that all $a_i \ge 1$.)

\paragraph{3.} 
An important observation is that \eqref{Ak axes} 
holds for {\em arbitrary} symmetric convex body $K$, 
in which case $a_i \sqrt{n}$ denote the semiaxes of an M-ellipsoid of $K$.
The M-ellipsoid is a deep concept in the modern convex geometry; 
it nicely reflects volumetric properties of convex bodies. 
For every symmetric convex body $K$ in $\R^n$ there exists 
an ellipsoid $\EE$ such that $|K| = |\EE|$ and $K$ can be covered
by at most $\exp(C_0 n)$ translates of $\EE$. Such ellipsoid $\EE$
is called an M-ellipsoid of $K$ (with parameter $C_0$). 
Its existence (with the parameter equal to an absolute constant)
was proved by V.Milman \cite{Mi 86}; for numerous consequences
see \cite{Pi}, \cite{Mi 98}, \cite{GM}.

\begin{fact} 
  Let $K$ be a symmetric convex body in $\R^n$ and $\EE$ be its
  M-ellipsoid with parameter $C_0$. Then 
  $$
  A_k(K) \le (CC_0)^{n/k} \Big( \prod_{i=1}^k a_i \Big)^{1/k},
  $$
  where $a_i \sqrt{n}$ are the semiaxes of $\EE$ in a nondecreasing order.
  In other words, $a_i$ are the singular values of a linear operator
  that maps $B_2^n$ onto $\EE$.
\end{fact}

\proof
The fact that $\EE$ is an M-ellipsoid of $K$ implies by standard
covering arguments that $(CC_0)^n |K \cap E| \ge |\EE \cap E|$ 
for all subspaces $E$ in $\R^n$, see e.g. \cite[Fact 1.1(ii)]{MTJ}.
Since $|K| = |\EE|$, we have $A_k(K) \le (CC_0)^{n/k} A_k(\EE)$, 
which reduces the problem to the examples of ellipsoids 
discussed above. 
\endproof

\paragraph{4.} 
A quantity similar to $A_k(K)$ and which equals 
$(\prod_{i=l}^{l+k} a_i)^{1/k}$ for the ellipsoid with 
nonincreasing semiaxes $a_i$
plays a central role in the recent work of Mankiewicz 
and Tomczak-Jaegermann \cite{MTJ}. They proved a volume ratio-type 
result for this quantity (for random non-coordinate subspaces $E$ 
in the Grassmanian) which works for $\dim E \le n/2$.

\paragraph{5.} 
Theorem \ref{cvr} follows from its more general dual counterpart
that allows to compute the combinatorial dimension of a set
in terms of its volume. 

Let $K$ be a set in $\R^n$. For $0 < k < n$, define
$$
a_k(K) = \min \left( \frac{|cK|}{|P_E K|} \right)^{1/\codim E}
$$
where the minimum is over the coordinate subspaces $E$, $\codim E \ge k$,
and $c>0$ is an absolute constant whose value will be discussed later.

\begin{theorem}                     \label{dual cvr}
  Let $K$ be a convex set in $R^n$. 
  Then for every integer $0 < k < n$, 
  $$
  v(K, a_k(K)) \ge n-k.
  $$
\end{theorem}

\proof
By applying an arbitrarily small perturbation to $K$, we can 
assume that the function $R \mapsto v(RK,1)$ maps $\R_+$ onto
$\{0, 1, \ldots, n\}$. Let $R$ be a solution to the equation 
$$
v(RK,1) = n - k.
$$
By Corollary \ref{vol cell content},
\begin{equation}                    \label{maxp}
  \big| \frac{1}{2}RK \big| 
  \le 1 + \max_P \# (\text{integer cells in $P(RK)$}) 
\end{equation}
where the maximum is over all coordinate projections $P$ in $\R^n$.
Since $v(RK,1) \ge 1$, the maximum in \eqref{maxp} is at least $1$. 
Hence there exists a coordinate projection $P = P_E$ onto a 
coordinate subspace $E$ such that 
$$
\big| \frac{1}{2}RK \big| 
\le 2 \# (\text{integer cells in $P_E(RK)$}).
$$
Since the number of integer cells in a set is bounded by its volume, 
$$
R^n \big| \frac{1}{2}K \big|
\le 2 |P_E(RK)| 
\le 2 R^{n-l} |P_E K|
$$
where $n-l = \dim E$.
It follows that 
$$
\frac{1}{R} \ge \left( \frac{|\frac{1}{4}K|}{|P_E K|} \right)^{1/l}
\ \ \ \text{and} \ \ \ 
v(K, \frac{1}{R}) = n-k.
$$
It only remains to note that by the maximal property of the 
combinatorial dimension, $n-l = \dim E \le n-k$; 
thus $l = \codim E \ge k$.
\endproof

\begin{lemma}                       \label{akak}
  For every integer $0 < k < n$, we have $A_k(K) \, a_k(nK^\circ) \ge 1$.
\end{lemma}

\proof
Let $L = nK^\circ$.
Fix numbers $0 < k \le l < n$ and a coordinate subspace $E$, $\codim E = l$.
Santalo and the reverse Santalo inequalities 
(the latter due to Bourgain and Milman \cite{BM}, see \cite{Pi} \S 7) 
imply that 
\begin{gather*}
  |L| \ge c_1^n |K|^{-1}, \\
  |P_E L| \le \left( \frac{C_1}{n-l} \right)^{n-l} |L^\circ \cap E|^{-1}
   = \left( \frac{C_1n}{n-l} \right)^{n-l} |K \cap E|^{-1}.
\end{gather*}
Then 
$$
\left( \frac{|cL|}{|P_E L|} \right)^{1/l}
\ge \left[ (c_1c)^n \left(\frac{n-l}{C_1 n}\right)^{n-l} 
    \frac{|K \cap E|}{|K|}
    \right]^{1/l}
\ge \left( \frac{|K \cap E|}{|(C_2/c) K|} \right)^{1/l}.
$$
Now take the minimum over $l \ge k$ and over $E$ to see that 
$a_k(L) \ge A_k(K)^{-1}$ if we choose $C = C_2/c$.
\endproof

\medskip \remark Theorem \ref{dual cvr} holds for general sets $K$ 
(not necessarily convex) if in the definition of $a_k(K)$ 
one replaces $|P_E K|$ by $|\cconv P_E K|$. 
The proof above easily modifies.

\paragraph{Proof of Theorem \ref{cvr}.}
By Theorem \ref{dual cvr} and Lemma \ref{akak}, 
$$
v(K^\circ, (n A_k(K))^{-1})
= v(nK^\circ, A_k(K)^{-1}) 
\ge n-k.
$$
By the symmetry of $K$, this means that exists an orthogonal 
projection $P_E$ onto a coordinate subspace $E$, $\codim E = k$, 
such that
$$
P_E (K^\circ) \supset P_E 
   \big( (n A_k(K))^{-1} 
   [-{\textstyle \frac{1}{2}}, {\textstyle \frac{1}{2}}]^n \big).
$$
Dualizing, we obtain 
$$
K \cap E \subset 2 A_k(K) B_1^n.
$$
The constant $2$ can be removed by increasing the 
value of the absolute constant $C$ in the definition of $A_k(K)$.
\endproof

\section{Volumes of the sets in the $L_p$ balls}        \label{s:Lp}

The classical Volume Ratio Theorem stated in the introduction 
is sharp up to an absolute sonstant $C$ (see e.g. \cite{Pi} \S 6). 
However, if we look at the factor $|CK|^{1/k} = C^{n/k} |K|^{1/k}$ 
which also appears in Corollary \ref{cvr cube inside intro}, 
then it becomes questionable 
whether the exponential dependence of the proportion $n/k$ 
is the right one. We will improve it in the dual setting 
to a {\em linear} dependence. The main result of this section computes the 
combinatorial dimension of a set $K$ (not even convex) in $\R^n$
in terms of its volume restricted to $B_p^n$. 
In other words, we are looking at the probability measure defined as
$$
\mu_p(K) = \frac{|K \cap B_p^n|}{|B_p^n|}.
$$

\begin{theorem}                     \label{mup}
  Let $K$ be a set in $\R^n$ and
  $1 \le p \le \infty$. Then for every integer $0 < k \le n$
  one has 
  \begin{equation}              \label{v vol}
  v(K,t) \ge n-k
  \ \ \ \text{for}\ \ \ 
  t = c \Big(\frac{k}{n}\Big) \mu_p(K)^{1/k}.
  \end{equation}
\end{theorem} 

\remarks 
{\bf 1.} The result is sharp up to an absolute constant $c$.
An appropriate example will be given after the proof.

{\bf 2.} Corollary \ref{cvr cube inside} is an immediate consequence
of Theorem \ref{mup} by duality.

{\bf 3.} To compare Theorem \ref{mup} to the classical Volume Ratio Theorem, 
one can read \eqref{v vol} for convex bodies as follows:
$$
(*) \ \ 
\begin{array}{l}
  \text{There exists a coordinate projection $P$ of rank $n-k$ 
        so that $PK$}\\
  \text{contains a translate of the cube $P(tB_\infty^n)$ with
         $t = c (\frac{k}{n}) \mu_p(K)^{1/k}$,}
\end{array}
$$
while the classical Volume Ratio Theorem states that 
$$
(**) \ \ 
\begin{array}{l}
  \text{There is a random orthogonal projection $P$ of rank $n-k$ 
        so that $PK$}\\
  \text{contains a translate of the ball $P(tB_2^n)$ with
         $t = c^{n/k} \mu_2(K)^{1/k}$.}
\end{array}
$$
Beside the central fact of the existence of a {\em coordinate} 
projection in $(*)$, note also the linear dependence on the 
proportion $k/n$ (in contrast to the exponential dependence in $(**)$), 
and also the arbitrary $p$.

\qquad

For the proof of Theorem \ref{mup}, we will need to know that 
the volumes $w_p(n) = |B_p^n|$ approximately increase in $n$.

\begin{lemma}                   \label{vol monotone}
  $w_p(k) \le C w_p(n)$ provided $k \le n$.
\end{lemma}

\proof
We have 
$$
w_p(k) = k^{k/p} 
  \frac{(2 \Gamma(1+\frac{1}{p}))^k}{\Gamma(1+\frac{k}{p})},
$$
see \cite{Pi} (1.17).
Note that 
$$
a^{1/p} := 2 \Gamma(1+\frac{1}{p}) 
\ge 2 \min_{x > 0} \Gamma(x) \ge 1.76.
$$ 
We then use Stirling's formula 
$$
\Gamma(1+z) \approx e^{-z} z^{z+1/2}
$$
where $a \approx b$ means $ca \le b \le Cb$ for some 
absolute constants $c,C>0$.

Consider two cases.

\noindent 1. $k \ge p$. We have
\begin{equation}                \label{k>p}
  w_p(k) = \frac{(ak)^{k/p}}{\Gamma(1+\frac{k}{p})}
  \approx (ak)^{k/p} e^{k/p} 
    \Big(\frac{k}{p}\Big)^{-\frac{k}{p}-\frac{1}{2}}
  \approx (eap)^{k/p} \sqrt{\frac{p}{k}}.
\end{equation}
2. $k \le p$. In this case $\Gamma(1+\frac{k}{p}) \approx 1$, 
thus 
\begin{equation}                \label{k<p}
  w_p(k) \approx (ak)^{k/p}.
\end{equation} 

To complete the proof, we consider three possible cases.

\noindent (a) $k \le n \le p$. Here the lemma is trivially 
true by \eqref{k<p}.

\noindent (b) $k \le p \le n$. Here
\begin{align*}
  \frac{w_p(n)}{w_p(k)} 
  &\gtrsim \frac{(eap)^{n/p}}{(ak)^{k/p}} \sqrt{\frac{p}{n}}
   \ge a^{\frac{n-k}{p}} \sqrt{\frac{k}{n}}
   \ \ \ \text{(because $p \ge k$)} \\
  &\ge (1.76)^{n-k} \sqrt{\frac{k}{n}} \ge c > 0.
\end{align*}

\noindent (c) $p \le k \le n$. Here
$$
\frac{w_p(n)}{w_p(k)}
  \gtrsim  (eap)^{\frac{n-k}{p}} \sqrt{\frac{n}{k}}.
$$
Since $ep>1$, one finishes the proof as in case (b). 
\endproof

\qquad

\noindent {\bf Proof of Theorem \ref{mup}.}
We can assume that $K \subseteq B_p^n$. Let 
$$
u^n = \frac{|K|}{|B_p^n|}.
$$
By applying an arbitrarily small perturbation of $K$ we can 
assume that the function $R \mapsto v(RK,1)$ maps $\R_+$ onto
$\{0, 1, \ldots, n\}$. Then there exists a solution $R$ to the equation
$$
v(RK, 1) = n-k.
$$ 
The geometric results of the previous sections, such as 
Corollary \ref{vol cell content}
and Theorem \ref{counting cells}, contain absolute constant factors 
which would destroy the linear dependence on $k/n$.
So we have to be more careful and apply \eqref{average content}
together with \eqref{cconv conv} instead:
\begin{equation}                    \label{maximum}
|RK| \le 1 + \max_{x \in (0,1)^n} \sum_P 
         \# \big( \text{integer cells in $Px + \cconv P(RK)$} \big). 
\end{equation}
Since $v(RK,1) > 0$, there exists a coordinate projection $P$
such that 
$$
\max_{x \in (0,1)^n}
\max_P \# \big( \text{integer cells in $Px + \cconv P(RK)$} \big) \ge 1.
$$
Hence the maximum in \eqref{maximum} is bounded below by $1$
(for $x = 0$). Thus  
\begin{align*}
|RK| 
&\le 2 \max_{x \in (0,1)^n} 
       \sum_P \# \big( \text{integer cells in $Px + \cconv P(RK)$} 
         \big) \\
&\le 2 \max_{x \in (0,1)^n}
       \sum_{d=1}^{n-k} \sum_{\rank P = d} 
       \# \big( \text{integer cells in $Px + \cconv P(RK)$} \big) \\
&\le 2 \sum_{d=1}^{n-k} \sum_{\rank P = d} |\cconv P(RK)|
\end{align*}
because the number of integer cells in a set is bounded by its volume. 
Note that
$\cconv P(RK) \subset \conv P(RK) \subset R P (B_p^n)$ by the assumption. 
Then denoting by $P_d$ the orthogonal projection in $R^n$ onto $\R^d$, 
we have
$$
|RK| \le 2 \sum_{d=1}^{n-k} \binom{n}{d} R^d |P_d B_p^n|.
$$
Now note that $P_d B_p^n = (n/d)^{1/p} B_p^d$. Hence
\begin{equation}                    \label{RK}
|RK| \le 
2 \sum_{d=1}^{n-k} \binom{n}{d} 
                   \Big(\frac{n}{d}\Big)^{d/p} R^d w_p(d).
\end{equation}
Now $|RK| = R^n |K| = R^n u^n w_p(n)$ in the left hand side
of \eqref{RK} 
and $w_p(d) \le C w_p(n)$ in the right hand side of \eqref{RK}
by Lemma \ref{vol monotone}. After dividing \eqref{RK} through 
by $R^n w_p(n)$ we get
\begin{equation}                    \label{vn bounded}
u^n \le 2C \sum_{d=1}^{n-k} \binom{n}{d} \Big(\frac{n}{d}\Big)^{d/p} 
        R^{d-n}.
\end{equation}

Let $0 < \e < 1$.
There exists a $1 \le d \le n-k$ such that 
$$
\Big(\frac{n}{d}\Big)^{d/p} R^{d-n} 
\ge (2C)^{-1} \e^{n-d} (1-\e)^d u^n;
$$
otherwise \eqref{vn bounded} would fail by the Binomial Theorem. 
From this we get
$$
R \le (2C)^{\frac{1}{n-d}} \Big(\frac{n}{d}\Big)^{\frac{d}{p(n-d)}}
      \frac{1}{\e} \Big(\frac{1}{1-\e}\Big)^{\frac{d}{n-d}}
      u^{-\frac{n}{n-d}}.
$$
Define $\d$ by the equation $d = (1-\d)n$. We have
$$
R \le (2C)^{1/\d n} \big[ (1-\d)^{1/p}(1-\e) \big]^{-(\frac{1-\d}{\d})}
      \frac{1}{\e} u^{-1/\d}.
$$
Now we use this with $\e$ defined by the equation $n-k = (1-\e)n$. 
Since $d \le n-k$, we have $\e \le \d$, so 
$$
\big[ (1-\d)^{1/p}(1-\e) \big]^{-(\frac{1-\d}{\d})}
\le (1-\e)^\frac{2(1-\e)}{\e} < C  
\ \ \ \text{for $0 < \e < 1$.}
$$
Thus
$$
R \le \frac{C}{\e} u^{-1/\e}.
$$
Then for $t := C^{-1} \e u^{1/\e} \le \frac{1}{R}$ we have 
$v(K,t) \ge v(K, \frac{1}{R}) = n-k$.
\endproof

\begin{example}
For every integer $n/2 \le k < n$ there exists a 
coordinate convex body $K$ in $\R^n$ of arbitrarily small volume 
and such that for all $1 \le p \le n$
$$
v(K,t) > n-k  \ \ \ \text{implies} \ \ \ 
t < C \Big(\frac{k}{n}\Big) \mu_p(K)^{1/k}.
$$
\end{example}

\proof
Fix an $\e > 0$ and let $K$ be the set of all points 
$x \in B_p^n$ such that one has $|x(i)| \le \e$
for at least $k$ coordinates $i \in \{1,\ldots, n\}$. 
Then $K$ contains $\binom{n}{k}$ disjoint sets $K_A$ 
indexed by $A \subset \{1,\ldots,n\}$, $|A|=k$,
$$
K_A = \{ x \in B_p^n : \; \text{one has $|x(i)| \le \e$ iff $i \in A$} \}.
$$
For each $A$, write 
$$
K_A = ([-\e,\e]^A \times (\e I)^{A^c}) \cap B_p^n
$$
where $I = (-\infty, -1) \cup (1,\infty)$.
In the next line we use notation $f(\e) \asymp g(\e)$ if 
$f(\e)/g(\e) \to 1$ as $\e \to 0$ uniformly over $p \in [1,\infty]$.
We have
\begin{align*}
|K_A| 
&\asymp | ([-\e,\e]^A \times \R^{A^c}) \cap B_p^n|
 \asymp |[-\e,\e]^A| \times |B_p^n \cap \R^{A^c}| \\
&= (2\e)^k \Big| \Big(\frac{n}{n-k}\Big)^{1/p} B_p^{n-k} \Big|
 \ge (2\e)^k |B_p^{n-k}|.
\end{align*}
Thus there exists an $\e = \e(n,k) > 0$ so that 
$$
\mu_p(K) = \binom{n}{k} \mu_p(K_A) 
  \ge \binom{n}{k} (c\e)^k \frac{|B_p^{n-k}|}{|B_p^n|}.
$$
Now we need now to bound below the ratio of the volumes. 

CLAIM. $\frac{w_p(n-k)}{w_p(n)} \ge c^k$.

Consider two possible cases: 

(a) $p \ge n-k$. In this case $n/2 \le n-k \le p \le n$, and  
by \eqref{k>p} and \eqref{k<p} we have
\begin{align*}
\frac{w_p(n-k)}{w_p(n)} 
&= \frac{(a(n-k))^{\frac{n-k}{p}}}{(eap)^{\frac{n}{p}}}
  \sqrt{\frac{n}{p}} \\
&\ge \Big( \frac{n-k}{ean} \Big)^{\frac{n}{p}}
     \ \ \ \text{(since $p \le n$)} \\
&\ge \Big( \frac{1}{2ea} \Big)^2 
     \ \ \ \text{(since $p \ge n-k \ge n/2$)}
\end{align*}
which proves the claim in this case. 

(b) $p \le n-k$. Here 
\begin{align*}
\frac{w_p(n-k)}{w_p(n)} 
&= \frac{(eap)^{\frac{n-k}{p}}}{(eap)^{\frac{n}{p}}}
  \sqrt{\frac{n}{n-k}} \\
&\ge \frac{1}{2} (eap)^{-\frac{k}{p}}
     \ \ \ \text{(since $n/2 \le k \le n$)} \\
&\ge c^k.
\end{align*}
This proves the claim.

We have thus shown that $\mu_p(K) \ge \binom{n}{k} (c \e)^k$, so 
$$
\mu_p(K)^{1/k}  >  c \Big(\frac{n}{k}\Big) \e.
$$
On the other hand, no coordinate projection $PK$ 
of dimension exceeding $n-k$ can contain a translate of 
the cube $P[-t,t]^n$ for $t > \e$. Thus 
$$ 
v(K,t) > n-k  \ \ \ \text{implies} \ \ \ 
t \le \e < C \Big(\frac{k}{n}\Big) \mu_p(K)^{1/k}.
$$
Note also that the volume of $K$ can be made 
arbitrarily small by decreasing $\e$.
\endproof

The same example also works for $p=\infty$.

\section{Duality for diameters of coordinate sections}      

Here we prove Theorem \ref{coordinate}. Formally, 
$$
r_k(K) =  \frac{2}{\sqrt{n}} \min_{|I|=k} \max_{x \in K}
  \sum_{i \in I} |x(i)|.
$$

\begin{theorem}
  Let $K$ be a symmetric convex body in $\R^n$.
  For any $\e > 0$ and for any two positive integers $k$ and $m$ 
  satisfying $k + m \le (1-\e)n$
  one has 
  $$
  r_k(K) \, r_m(K^\circ) \le C^{1/\e}.
  $$
\end{theorem}
\noindent The proof is based on Corollary \ref{cvr cube inside}.

\qquad

\proof
Define $\d$ and $\l$ as folows: $k = (1-\d)n$, $m=(1-\l)n$. 
Then $\d + \l - 1 > \e$. 
Let $t_1, t_2 > 0$ be parameters, and define
$$
K_1  =  \conv \Big( K \cup t_1 n^{-1/2}B_\infty^n \Big) 
        \cap \frac{1}{t_2} n^{-1/2}B_1^n.
$$
Consider two possible cases:

(1) $|K_1| \le |n^{-1/2}B_\infty^n|$. 
Since $K_1$ contains $t_1 n^{-1/2}B_\infty^n$, we have 
$$
\frac{\sqrt{n}}{t_1} K_1 \supset B_\infty^n
\ \ \ \text{and} \ \ \
\Big| \frac{\sqrt{n}}{t_1} K_1 \Big|  
\le  \Big| \frac{1}{t_1} B_\infty^n \Big|
= \Big(\frac{2}{t_1}\Big)^n.
$$
Corollary \ref{cvr cube inside}
implies the existence of a subspace $E$, $\dim E = (1 - \d) n$, 
such that 
$$
\frac{\sqrt{n}}{t_1} K_1 \cap E  
\subset \Big( \frac{C}{t_1} \Big)^{1/\d} B_1^n.
$$
Multiplying through by $t_1 / \sqrt{n}$ and recalling the definition
of  $K_1$, we conclude that 
\begin{equation}                                           \label{1}
K \cap E \cap \frac{1}{t_2} n^{-1/2}B_1^n
  \subset  t_1 \Big( \frac{C}{t_1} \Big)^{1/\d} n^{-1/2}B_1^n.
\end{equation}

(2) $|K_1| > |n^{-1/2}B_\infty^n|$. 
Note that
$$
K_1^\circ  
=  \conv \Big[ \Big( K^\circ \cap \frac{1}{t_1} n^{-1/2}B_1^n \Big) 
         \cup t_2 n^{-1/2}B_\infty^n \Big].
$$
By Santalo and reverse Santalo inequalities, 
$$
|K_1^\circ|  <  |C n^{-1/2}B_1^n|.  
$$
Since $K_1^\circ$ contains $t_2 n^{-1/2}B_\infty^n$, we have 
$$
\frac{\sqrt{n}}{t_2} K_1^\circ \supset B_\infty^n
\ \ \ \text{and} \ \ \
\Big| \frac{\sqrt{n}}{t_2} K_1^\circ \Big|  
\le  \Big| \frac{C}{t_2} B_1^n \Big|
\le \Big(\frac{C_2}{t_2}\Big)^n.
$$
Arguing similarly to case (1) for $K^\circ$, we find
a subspace $F$, $\dim F = (1 - \l) n$, and such that 
\begin{equation}                                           \label{2}
K^\circ \cap F \cap \frac{1}{t_1} n^{-1/2}B_1^n
  \subseteq  t_2 \Big( \frac{C}{t_2} \Big)^{1/\l} n^{-1/2}B_1^n.
\end{equation}
Looking at \eqref{1} and \eqref{2}, we see that 
our choice of $t_1, t_2$ should be so that 
$$
t_1 \Big( \frac{C}{t_1} \Big)^{1/\d} = \frac{1}{2t_2}, \ \ \ \
t_2 \Big( \frac{C}{t_2} \Big)^{1/\l} = \frac{1}{2t_1}.
$$
Solving this for $t_1$ and $t_2$ we get
$$
\frac{1}{2t_1} = \frac{1}{\sqrt{2}} C^{\frac{\d-\l+1}{\d+\l-1}} =: R_1, \ \ \ \
\frac{1}{2t_2} = \frac{1}{\sqrt{2}} C^{\frac{\l-\d+1}{\d+\l-1}} =: R_2.
$$
Then \eqref{1} becomes 
$$
K \cap E   \subseteq R_2 n^{-1/2}B_1^n
$$
and \eqref{2} becomes 
$$
K^\circ \cap F   \subseteq R_1 n^{-1/2}B_1^n.
$$
It remains to note that 
$$
R_1 R_2 = \frac{1}{2} C^{2/(\d+\l-1)} < \frac{1}{2} C^{2/\e}.
$$
This completes the proof.
\endproof

\end{document}